# ON THE MAXIMAL SECTORIAL OPERATORS

Z.I. Ismailov and R. Öztürk

Karadeniz Technical University, Faculty of Arts and Sciences
Department of Mathematics
61080 Trabzon, Turkey

e-mail: zameddin@yahoo.com

**Abstract**

In this work, firstly the maximal sectorial linear relations are described. Later on, the discreteness of the spectrum of the linear maximal sectorial operators and asymptotical behaviour of the eigenvalues of such operators in terms of the eigenvalues of its real part are investigated. Finally, it is obtained one result for the differential operators for first order in the Hilbert space of vector functions in finite interval.

## 1. Introduction

The general theory of linear selfadjoint and normal operators, and their spectral analysis in a Hilbert space has been studying by many mathematicians (for example, see [1-4]).

Some problems for the linear unbounded selfadjoint operators, such as accretivity, sectoriality, the structure and discreteness of spectrum, the asymtotical behavior of eigenvalues at infinity, completeness of eigenvectors and so on, are very important for physical and technical problems (for example, see [5-7] and references in it).

Numerical range and sectoriality of operators have been investigated by many authors (for example, see [8-10]) and a full account of the basic functional calculus theory of these operators, applications of this theory to the theory of differential operators and their spectral theory have been given in work [10].

In addition, there are many applications of this theory to the different models (see [11], [12] and references in it).

In this paper, $H$ stands for a separable Hilbert space, $A_R$ and $A_I$ represent real and imaginary parts of a linear closed operator $A: D(A) \subseteq H \to H$, respectively, $B(H)$ refers to the family of linear bounded operators, $C_p(H)$, $1 \leq p < \infty$ is the Schatten-von Neumann class of compact operators, $C_\infty(H)$ stands for the set of compact operators in $H$, $\sigma(A)$ and $\rho(A)$ refer to spectrum and the resolvent sets of an operator $A: D(A) \subseteq H \to H$ respectively.

## 2. Description of m-sectorial operators

**Definition 2.1:** Let $H$ be a Hilbert space. The linear subspace
$$\theta = \{\{x, x'\} : x, x' \in H\} \subset H \oplus H$$
is called a linear m-sectorial linear relation if the following conditions are satisfied:
(1) For every $\{x, x'\}$ is true $\operatorname{Re}(x', x)_H \geq 0$, i.e. $\theta$ is a linear accretive relation;
(2) $\theta$ is a m-accretive relation;
(3) For any $\varphi \in \left[0, \dfrac{\pi}{2}\right]$ holds
$$\cos \varphi |\operatorname{Im}(x', x)_H| \leq \sin \varphi \operatorname{Re}(x', x)_H;$$

In this case for the m-sectorial linear relation there is a semi-angle $\varphi \in \left[0, \dfrac{\pi}{2}\right]$ such that
$$W(\theta) = \{(x', x)_H \in \mathbb{C} : \{x, x'\} \in \theta\} \subseteq \{z \in \mathbb{C} : |\arg z| \leq \varphi\}.$$

It is easy to prove that is true the following assertion.

**Theorem 2.2:** A linear relation $\theta \subset H \oplus H$ is a m-sectorial relation if and only if are satisfy the following properties:
(1) $\theta$ is a m-accretive relation;
(2) $\theta_\varphi := \{\{x, e^{i\varphi} x'\} : \{x, x'\} \in \theta\}$ is a m-dissipative relation;
(3) $\theta_{-\varphi} := \{\{x, e^{-i\varphi} x'\} : \{x, x'\} \in \theta\}$ is a m-accumulative relation.

Now can be describe the m-sectorial linear relation.

**Theorem 2.3:** If $\theta = \{\{x, x'\} : x, x' \in H\} \subset H \oplus H$,
$$W(\theta) \subset \{z \in \mathbb{C} : |\arg z| \leq \varphi\}, \ \varphi \in \left[0, \dfrac{\pi}{2}\right]$$
is a m-sectorial linear relation in $H$, then there are linear contraction operators $K, V, W$ in $H$, such that are satisfied the conditions
$$(K - E) x' + (K + E) x = 0,$$
$$(V - E) e^{i\varphi} x' + i(V + E) x = 0,$$
$$(W - E) e^{-i\varphi} x' - i(W + E) x = 0$$
and vice versa. Here the linear operators $K, V, W$ in $H$ are determined by the linear relation $\theta$ uniqely, i.e.
$$\theta = \theta(K, V, W).$$

**Proof:** Since $\theta$ is a maximal sectorial linear relation, the by the Theorem 2.2 $\theta, \theta_\varphi$ and $\theta_{-\varphi}$ are m-accretive, m-dissipative and m-accumulative relation in $H$ respectively. In this situation the linear relations

$$\theta = \theta(\varphi) \subset H \oplus H,$$

$$\theta_\varphi = \{\{x, e^{i\varphi} x'\} : \{x, x'\} \in H\} \text{ and}$$

$$\theta_{-\varphi} = \{\{x, e^{-i\varphi} x'\} : \{x, x'\} \in H\}$$

are described with above equations (see [6]). Uniqueness of the contractions $K, V, W$ it is clear.

On the contrary, it is clear that linear relation which is determined with solutions of the above writing equations is m-sectorial.

According to the m-sectorial linear relation can be defined a m-sectorial linear operators in any Hilbert space.

**Definition 2.4:** A densely defined linear closed operator $T$ in Hilbert space $H$ is called a m-sectorial, if:

(1) $T$ is a m-accretive operator in $H$,

(2) $W(T) = \{(Tx, x) \in \mathbb{C} : x \in D(T)\} \subset \{z \in \mathbb{C} : |\arg z| \leq \varphi\}$, $\varphi \in \left[0, \dfrac{\pi}{2}\right]$

(or $\cos\varphi |\text{Im}(Tx, x)| \leq \sin\varphi \text{Re}(Tx, x)$, $\varphi \in \left[0, \dfrac{\pi}{2}\right]$ ).

In this case we will write $T \in S_\varphi(H)$, $\varphi \in \left[0, \dfrac{\pi}{2}\right]$.

It is easy to prove that $T \in S_\varphi(H)$ if and only if $T$ m-accretive, $e^{i\varphi}T$ is m-dissipative and $e^{-i\varphi}T$ is a m-accumulative operators in $H$.

Now it is investigated the discreteness of the spectrum of the linear m-sectorial operators. In fact is true the following cleam.

**Theorem 2.5:** If $T \in S_\varphi(H)$ and for any $\alpha > 0$

$$(T_R + \alpha)^{-1} \in C_p(H), \ p \geq 1,$$

then

$$(T + \alpha)^{-1} \in C_p(H)$$

and true the asymptotical formula for the eigenvalues of the operator $T$

$$\lambda_n(T) \sim \text{Re}\,\lambda_n(T) \text{ as } n \to \infty.$$

**Proof:** Indeed, since the operator $T$ is m-accretive, the for any $\alpha > 0$, $(T_R + \alpha)^{-1} \in B(H)$. On the other hand can be written

$$T + \alpha = (T_R + \alpha)\left[(T_R + \alpha)^{-1} + i(T_R + \alpha)^{-1} T_I (T_R + \alpha)^{-1}\right](T_R + \alpha). \quad (2.1)$$

In this case for the numerical range $W(P)$ of the operator

$$P = (T_R + \alpha)^{-1} + i(T_R + \alpha)^{-1} T_I (T_R + \alpha)^{-1}$$

it is obtained that

$$(Px, x) = \left((T_R + \alpha)^{-1} x + i(T_R + \alpha)^{-1} T_I (T_R + \alpha)^{-1} x, x\right)$$

$$= \left((T_R + \alpha)^{-1} x, (T_R + \alpha)\left((T_R + \alpha)^{-1}\right)x\right) + i\left(T_I \left((T_R + \alpha)^{-1} x\right), (T_R + \alpha)^{-1} x\right),$$

$x \in D(P)$. From this relation
$$\operatorname{Re}(Px,x)_H = (y,(T_R+\alpha)y) = \alpha(y,y)+(T_R y,y)$$
$$\geq \alpha(y,y),$$
$y = (T_R+\alpha)^{-1} x$, $x \in D(P)$, $\alpha > 0$, i.e.
$$W(P) \subset \left\{ z \in \mathbb{C} : |\arg(z-\alpha)| \leq \varphi \leq \frac{\pi}{2} \right\}.$$
Consequently,
$$\left[ (T_R+\alpha)^{-1} + (T_R+\alpha)^{-1} T_I (T_R+\alpha)^{-1} \right]^{-1} \in B(H).$$
From this relation, equality (2.1) and property $BC_p \subset C_p$ it is established that $(T+\alpha)^{-1}$ exist and belong to $C_p(H)$, $p \geq 1$.

Now let us a operator $T \in S_\varphi(H)$ and $T$ have discrete spectrum in $H$. Assume that
$$\lambda \in \sigma_p(T) \text{ and } \lambda = \lambda^r + i\lambda^i, \ \lambda^r, \lambda^i \in \mathbb{R}.$$
In this case, there is a any eigenvector $x_\lambda \in H$, $x_\lambda \neq 0$, such that $Tx_\lambda = \lambda x_\lambda$. From this it is obtained
$$(Tx_\lambda, x_\lambda) = \lambda^r(x_\lambda, x_\lambda) + i\lambda^i(x_\lambda, x_\lambda)$$
and consequently
$$\lambda^r = \operatorname{Re}\frac{(Tx_\lambda, x_\lambda)}{(x_\lambda, x_\lambda)}, \ \lambda^i = \operatorname{Im}\frac{(Tx_\lambda, x_\lambda)}{(x_\lambda, x_\lambda)}.$$
Furthermore, it is easy to see that
$$|\lambda_n^i|^2 \leq |\lambda_n|^2 = |\lambda_n^r|^2 + |\lambda_n^i|^2 \leq (1+\tan^2 \varphi)|\lambda_n^r|^2, \ n \geq 1.$$
On the other words, it is true the following asymptotics
$$\lambda_n(T) \sim \lambda_n^r(T) = \operatorname{Re}\frac{(Tx_{\lambda_n}, x_{\lambda_n})}{(x_{\lambda_n}, x_{\lambda_n})} \text{ as } n \to \infty.$$
In special case, from this theorem implies the following result.

**Corollary 2.6:** If $T \in S_\varphi(H)$, $T$ linear normal opeator in $H$ and for some $\alpha > 0$ $(T_R+\alpha)^{-1} \in C_p(H)$, $p \geq 1$, then $(T+\alpha)^{-1} \in C_p(H)$, $p \geq 1$ and the following asimptotics is true
$$\lambda_n(T) \sim \lambda_n(T_R) \text{ as } n \to \infty.$$

## 3. Differential Operators of First Order

In this section will be study a m-sectoriality of extension of the minimal operator generated by linear differential operator expression in form
$$l(u) = u' + Au, \qquad (3.1)$$
where $A \in S_\varphi(H)$, $0 \leq \varphi < \frac{\pi}{2}$, in the Hilbert space $L_2(H,(a,b))$ of vector functions in finite interval.

It is clear that formally adjoint expression of (3.1) in the Hilbert space $L_2(H,(a,b))$ is in form
$$l^+(v) = -v' + A^* v. \qquad (3.2)$$

Let us define $L_0'$ on the dense $L_2(H,(a,b))$ set of vector-functions $D_0'$,
$$D_0' := \left\{ u \in L_2(H,(a,b)) : u(t) = \sum_{k=1}^{n} \varphi_k(t) f_k, \; \varphi_k \in C_0^\infty(a,b), \; f_k \in D(A), \; k=1,2,\ldots,n, \; n \in \mathbb{N} \right\}$$
as $L_0' u = l(u)$. Since the operator $\operatorname{Re}(Af, f) \geq 0$, $f \in D(A)$, then operator $L_0'$ is accretive, that is $\operatorname{Re}(L_0' u, u)_{L_2(H,(a,b))} \geq 0$, $u \in D_0'$. Hence the operator $L_0'$ has a closure in $L_2(H,(a,b))$ [13]. The closure of $L_0'$ in $L_2(H,(a,b))$ is called the minimal operator generated by differential expression (3.1) and it is denoted by $L_0$.

In a similar way the minimal operator $L_0^+$ in $L_2(H,(a,b))$ can be construct generated by differential-operator expression (3.2). The adjoint operator of $L_0^+(L_0)$ in $L_2(H,(a,b))$ is called the maximal operator generated by (3.1) ((3.2)) and is denoted by $L$ ($L^+$).

Before of all it is investigated the relationship between the accretivity property of the minimal operator $L_0$ in $L_2(H,(a,b))$ and accretivity property of the operator coefficient $A$ in $H$.

**Theorem 3.1:** In order to the minimal operator $L_0$ in $L_2(H,(a,b))$ is accretive the necessary and sufficient condition is accretivity of the operator coefficient $A$.

**Proof:** Let us the operator $L_0$ is accretive in $L_2(H,(a,b))$. In this case the special functions $u(t) = \varphi(t) f$, $\varphi(t) \in \overset{0}{W_2^1}(a,b)$, $f \in D(A)$ belong to the $D(L_0)$ and for the such functions we have
$$L_0 u = u' + Au = \varphi'(t) f + \varphi(t) Af.$$
It is clear that
$$\|\varphi'(t) f\|_{L_2}^2 = \|\varphi'(t)\|_{L_2}^2 \|f\|_H^2 < \infty$$
and
$$\|\varphi(t) Af\|_{L_2}^2 = \|\varphi(t)\|_{L_2}^2 \|Af\|_H^2 < \infty.$$
On the other hand for the functions in form $u(t) = \varphi(t) f$, $\varphi(t) \in \overset{0}{W_2^1}(a,b)$, $\|\varphi\|_{L_2} \neq 0$, $f \in D(A)$ after the simple transformations we have
$$\operatorname{Re}(L_0 u, u)_{L_2} = \|\varphi\|_{L_2}^2 \operatorname{Re}(Af, f)_H \geq 0.$$
Therefore for any $f \in D(A)$ we have $\operatorname{Re}(Af, f) \geq 0$, i.e. $A$ is an accretive in $H$.

Now let $A$ be a accretive operator in $H$. In this case in the dense set of vector functions
$$D_0' := \left\{ u(t) = \sum_{k,n=1}^{m} \varphi_k(t) f_n : \varphi_k \in \overset{0}{W_2^1}(a,b), \; f_n \in D(A), \; k,n = 1,\ldots,m, \; m \in \mathbb{N} \right\}$$

it is obtained that
$$\operatorname{Re}(L_0 u, u)_{L_2} = \operatorname{Re}(Au, u)_{L_2}, \ u \in D_0'.$$
That is the restriction $L_0|_{D_0'}$ of the operator $L_0$ to the linear dense set $D_0'$ in $L_2(H,(a,b))$ is a accretive in $L_2(H,(a,b))$. From this and result in [13] it is obtained that the minimal operator $L_0$ is accretive in $L_2(H,(a,b))$.

**Theorem 3.2:** Let $A = A^* \geq 0$ and $AW_2^1(H,(a,b)) \subset L_2(H,(a,b))$. A minimal operator $L_0$ generated by linear differential-operator expression in form $l(u) = u' + Au$ in the Hilbert space $L_2(H,(a,b))$ has no any m-sectorial extension $\tilde{L} \in S_\varphi(L_2)$, $0 \leq \varphi < \dfrac{\pi}{2}$.

**Proof:** Assume that the minimal operator $L_0$ has an extension $\tilde{L}$ such that $\tilde{L} \in S_\varphi(L_2)$, $0 \leq \varphi < \dfrac{\pi}{2}$. In this case by the definition before of all $\tilde{L}$ is a m-accretive. Consequently, $\tilde{L}$ generates by $l(u)$ and boundary condition
$$u(a) = Ku(b),$$
where $K$ is a contraction operator in $H$ [14]. Now introduce the following functions in form
$$u_n(t) = \frac{1}{\sqrt{2(b-a)}} \left[ \exp\left(\frac{2n\pi i}{b-a}(t-a)\right) - 1 \right] f, \ f \in D(A), \ \|f\| = 1, \ Af \neq 0, \ n \geq 1.$$
It is clear that $u_n \in D(L_0) \subset D(\tilde{L})$, $\|u_n\| = 1$, $n \geq 1$. Simple calculations show that
$$(\tilde{L}u_n, u_n) = (Af, f) + \frac{\pi n}{b-a} i, \ n \geq 1.$$
Hence,
$$\frac{|\operatorname{Im}(\tilde{L}u_n, u_n)|}{\operatorname{Re}(\tilde{L}u_n, u_n)} = \frac{n\pi}{(b-a)\|Af\|^2} \to \infty \text{ as } n \to \infty.$$
This is a contraction. This completes a proof of the theorem.

**Remark 3.3:** From last theorem and work [7] implies that a minimal operator have not any normal m-sectorial extension.

# References

1. Von Neumann, J., Allgemeine eigenwerttheorie hermitescher funktionaloperatoren, Math. Ann., 102 (1929-1930) 49-131.

2. Birman, M. S., Solomjak, M. Z., Spectral theory of self-adjoint operators in Hilbert space, D. Reidel Publishing Company Dordrecht, Holland, 1987.

3. Dunford, N., Schwartz, J.T., Linear Operators, , v. II, Interscience Publisher, INC., New York, 1963.

4. Coddington, E.A., Extension theory of formally normal and symetric subspaces, Mem. Amer. Math. Soc., 134 (1973) 1-80.

5. Arlinskii, Yu. M., Boundary conditions for maximal sectorial extensions of sectorial operators, Math. Nachr., 209(2000) 5-36.

6. Gorbachuk, V.I., Gorbachuk, M.L., Boundary value problems for operator differential equations, Kluwer Academic, Dordrecht, 1991.

7. Ismailov, Z.I., Compact inverses of first-order normal differential operators, J. Math. Anal. Appl., 320 (2006) 266-278.

8. Kato, T., Perturbation theory of linear operators, Springer-Verlag, Berlin, Heidelberg, 1995.

9. Schechter, M., Principles of functional analysis, New York-London, Academic Press, 1971.

10. Haase, M., The functional calculus for sectorial operators, Birkhauser Verlag, Basel-Boston-Berlin, 2006.

11. Mihailets, V.A., On the solvable and sectorial boundary value problems for Sturm-Liouville operator equations, Ukr. Math. Journal, v.26, no.4 (1974) 450-459 (in Russian).

12. Malamud, M., Operator holes and extensions of sectorial operators and dual pairs of contractions, Math. Nach., 279 (2006) 625-655.

13. Sz.-Nagy, B., Foias, C., Harmonic analysis of operators on Hilbert space, Amsterdam-Budapest, North-Holland Publishing Co., 1970.

14. Levchuk, V.V., Smoothly maximal dissipative boundary value problems for parabolic equation in Hilbert space, Ukr. Math. Journal, v.36, no.4 (1983) 502-507 (in Russian).